\documentclass
{amsproc}
\usepackage{palatino, mathpazo}
\usepackage{enumitem, MnSymbol} 
\usepackage{hyperref}
\hypersetup{
	colorlinks=true,
	linkcolor=blue,
	citecolor=blue,
	urlcolor=blue
}
\newtheorem{theorem}{Theorem}
\newtheorem{proposition}[theorem]{Proposition}
\newtheorem{corollary}[theorem]{Corollary}
\theoremstyle{remark}
\newtheorem{remark}[theorem]{Remark}
\usepackage{geometry}
\newcommand{\overbar}[1]{\mkern 2mu\overline{\mkern-2mu#1\mkern-1.5mu}\mkern 1.5mu}
\newcommand{\p}{\partial}
\newcommand{\pbar}{\overline{\partial}}
\newcommand{\RR}{\mathbb{R}}
\newcommand{\CC}{\mathbb{C}}
\newcommand{\X}{\mathcal{X}}
\newcommand{\bbeta}{\overline{\beta}}
\numberwithin{equation}{section}

\makeatletter
\let\sv@thm\@thm
\def\@thm{\let\indent\relax\sv@thm}

\renewcommand\subsection{%
  \@startsection{subsection}%  name
    {2}% level
    {0em}% indent
    {2.4ex plus 0.1ex minus -0.1ex}% beforeskip
    {.08ex plus 0.05ex minus -0.05ex}% afterskip
    {\bf}% style
  }
\makeatother

\title[]{Plurisubharmonicity of the Dirichlet energy \\and deformations of polarized manifolds}
\author[]{\sc Che-Hung~Huang\vspace{-4ex}}
\address{
Department of Mathematics, Purdue University,
West Lafayette, IN 47907, USA}
\email{huan1160@purdue.edu}
\thanks{This research was partially supported by NSF grant DMS-1764167}

\begin{document}\vspace*{1.25cm}\maketitle

\bigskip\bigskip
\noindent\makebox[\textwidth][c]{
\begin{minipage}{0.835\textwidth}
 \footnotesize {\sc Abstract.} 
   We show that if $\{M_t\}_{t\in \Delta}$ is a polarized family of compact K\"ahler manifolds over the open unit disk $\Delta$, if $N$ is a Riemannian manifold of nonpositive complexified sectional curvature, and if $\{\phi_t:M_t\to N\}_{t\in \Delta}$ is a smooth family of pluriharmonic maps, then the Dirichlet energy $E(\phi_t)$ is a subharmonic function of $t\in\Delta$.
 \end{minipage} }
\smallskip\vspace{-1ex}

 \thispagestyle{empty}
 \section{Introduction}
  It is well known that harmonic maps from compact K\"ahler manifolds have the following rigidity property:
 \begin{theorem}[\text{Siu \cite{Siu}, Sampson \cite{Sampson}}]\label{theorem 1}Let $\phi:M\to N$ be a harmonic map from a compact K\"ahler manifold $(M,g)$ into a Riemannian manifold $(N,h)$. Suppose that $(N,h)$ satisfies the curvature condition: \begin{equation}\label{1.1}
  R^N(X,Y,\overline{X}, \overline{Y})\leq 0\ \ \, \text{for all $X,Y\in T_{\CC, p} N$ and all $p\in N$}. \end{equation} Then $\phi$ is pluriharmonic.
 \end{theorem}
 
   In this paper, we consider the effect of deforming the domain manifold $M$ and the map $\phi$, and we provide an application of a Bochner type identity to noncompact K\"ahler manifolds. 
 It turns out that Theorem \ref{theorem 1} is closely related to plurisubharmonicity of the Dirichlet energy. The statement of our main theorem is the following:
  \begin{theorem}\label{theorem 2}
  Suppose $(\pi:\X\to\Delta, \widetilde{\omega})$ is a polarized family of compact K\"ahler manifolds parametrized by the open unit disk $\Delta$. Suppose $M_t=\pi^{-1}(t)$ for every $t\in\Delta$, $(N, h)$ is as in Theorem \ref{theorem 1}, and $\{\phi_t: M_t\to N\}_{t\in\Delta}$ is a smooth family of pluriharmonic maps. Let $E:\Delta\to\RR$ be the function that assigns to each $t\in \Delta$ the energy $E(\phi_t)$ of $\phi_t$. Then $E$ is subharmonic.
  \end{theorem}
   Recall that a {\it polarized family of compact K\"ahler manifolds} parametrized by a reduced complex space $S$ is a pair $(\pi:\X\to S,\widetilde{\omega})$, where  $\pi:\X\to S$ is a proper smooth morphism of complex spaces, with connected fibers, and  $\widetilde{\omega}$ is an element of $\Gamma(S, R^2\pi_*\RR)$ that induces on each fiber $M_t=\pi^{-1}(t)$ a K\"ahler class $\widetilde{\omega}_t\in H^2(M_t,\RR)$.\\
   \indent Under the assumptions of Theorem \ref{theorem 2}, the energy $E(\phi_t)$ of $\phi_t$, defined with respect to the K\"ahler class $\widetilde{\omega}_t\in H^2(M_t,\RR)$ and the metric $h$, depends smoothly on $t\in\Delta$. We shall derive an explicit expression for $\p^2 E/\p t\p\overline{t}$ (formula \eqref{3.3}) which shows that $\p^2 E/\p t\p\overline{t}\geq 0$ on $\Delta$. 
   If the open unit disk $\Delta$ in the assumptions of Theorem \ref{theorem 2} is replaced by a domain $\Omega$ in $\CC^k$, with $k\geq 2$, it follows immediately that the function $E:\Omega\to \RR$ is plurisubharmonic.\\
 \indent The subject of plurisubharmonicity of the Dirichlet energy goes back at least as far as the\linebreak work of Tromba \cite{Tromba1, Tromba2}. The results of Tromba show that Theorem \ref{theorem 2} holds in the case where $M_t$ are compact Riemann surfaces of genus $g\geq 2$, $N=M_0$, and $\phi_0$ is the identity map on $M_0$. Toledo \cite{Toledo} extended Tromba's analysis to the case where $N$ is a Riemannian manifold satisfying $\eqref{1.1}$. Later, Kim et al.\ \cite{Kim} recovered Toledo's results using the approach of fiber integration.
  The case where $M_t$ are K\"ahler--Einstein manifolds of negative Ricci curvature was treated in \cite{Zhang}. Our method extends that of \cite{Kim}, and the proof of Theorem \ref{theorem 2} builds on the following observation: 
     First,
      the energy function $E$ 
     can locally be expressed as a fiber\linebreak integral. Since the 
     fiber integration $\pi_*$ commutes with both $\p$ and $\pbar$ (the family  $\pi:\X\to \Delta$ being holomorphic), this yields a fiber integral formula for $\sqrt{-1}\,\p\pbar E$. It is rather interesting to note that the last formula links the Levi form of $E$ to a Bochner type identity closely related to the rigidity theorem of Siu and Sampson. \bigskip
  
   \section{Preliminaries}
   In this section we review an argument for Theorem \ref{theorem 1} following the ideas in \cite{Siu} (see also \cite{Sampson}); such formulation will be meaningful for our interpretation of the Levi form of $E$, and\linebreak it will provide a clearer picture of the connection between Theorem \ref{theorem 1} and \ref{theorem 2}. 
  \\\indent   
  Let $\phi: M\to N$ be a smooth map from an $m$-dimensional complex manifold $M$ into a Riemannian manifold $(N, h=h_{ij}\,dy^i\otimes dy^j)$. The complexified differential $d\phi:T_\CC M\to \phi^*T_\CC N$ can be decomposed as a sum $d\phi=\p \phi+\pbar \phi$, with $\p\phi=\p\phi^i\otimes\p/\p y^i$ and $\pbar\phi=\pbar\phi^j\otimes\p/\p y^j$, according to the decomposition $T_\CC M
  = T_M\oplus \overbar{T_M}$, where $T_M$ is the holomorphic tangent bundle of $M$. Following \cite{Ohnita, Siu}, we define \[\varepsilon(\phi)= \sqrt{-1}\,h_{ij}\,\p \phi^i\wedge\overline{\p}\phi^j,\] which is a real $(1,1)$-form on $M$. Let $D$ be the exterior covariant derivative on $\phi^* T_{\CC} N$ induced by the Levi--Civita connection on $N$. Let $D'$ and $D''$ be the $(1,0)$- and $(0,1)$-components of $D$, respectively. Then \[D'\p \phi = 0 = D''\pbar \phi\ \ \ \  \text{and}\ \ \ \ D'\pbar \phi =\left(\p\pbar\phi^k+\Gamma^k_{ij}\,\p\phi^i\wedge\pbar\phi^j\right)\otimes\p/\p y^k= - D''\p \phi, \]where $\Gamma^k_{ij}$ are the Christoffel symbols of $(N,h)$. Since $D$ is compatible with the Hermitian structure on $\phi^*T_\CC N$ induced by $h$, we have \begin{equation}\label{2.1} \pbar\varepsilon(\phi)= \sqrt{-1}\,h_{ij}\,D''\p \phi^i\wedge\overline{\p}\phi^j, \end{equation}and hence \[\p\pbar\varepsilon(\phi)=\sqrt{-1}\,h_{ij}D'\,D''\p \phi^i\wedge\overline{\p}\phi^j+\sqrt{-1}\,h_{ij}\,D''\p \phi^i\wedge D'\overline{\p}\phi^j.\] We denote by ${R^l}_{kij}$ the components of the $(1,3)$-curvature tensor on $N$; in normal coordinates $(y^i)$ centered at a point $p\in N$, we have ${R^l}_{kij}=\p\Gamma^l_{jk}/\p y^ i-\p\Gamma^l_{ik}/\p y^j$ at $p$. Observe that  \[D'D''\p\phi=-D'D'\pbar\phi=-\frac{1}{2!}\,{R^l}_{kij}\,\p\phi^i\wedge\p \phi^j\wedge\pbar\phi^k\otimes{\p}/{\p y^l},\] since $(D')^2$ is the $(2,0)$-component of the curvature $2$-form $D^2$. We thus obtain the following Bochner type identity (cf.\ \cite [Proposition 2]{Siu}): \begin{equation}\label{2.2}\p\pbar\varepsilon(\phi)=\frac{\sqrt{-1}}{2}\,R_{ijkl}\,\p \phi^i\wedge\p \phi^j\wedge\pbar \phi^k\wedge\pbar \phi^l+\sqrt{-1}\, h_{ij}\,D''\p \phi^i\wedge D'\pbar \phi^j.\end{equation} \indent Let $g$ be a Hermitian metric on $M$, with fundamental $2$-from $\omega$. Suppose that $\{Z_\alpha\}$ is a local unitary frame of $T_M$ with respect to $g$. 
 A straightforward computation yields the following formula (\cite[(1.3)]{Ohnita}):
 \begin{align}\label{2.3}\sqrt{-1}\, \p\pbar\varepsilon(\phi)\wedge\frac{\omega^{m-2}}{(m-2)!}=\bigg\{&-\sum_{\alpha,\beta=1}^m R^N(\phi_*Z_\alpha ,\, \phi_*Z_\beta,\, \phi_*\overbar{Z_\alpha},\, \phi_*\overbar{Z_\beta})\\\nonumber &+|D''\p \phi|^2-\lvert\text{tr}_g D''\p\phi\rvert^2\bigg\}\,\frac{\omega^m}{m!},\end{align} where $\text{tr}_g D''\p\phi=\sum_{\alpha=1}^m (D_{\overbar{Z_\alpha}}\p \phi)(Z_\alpha)$.
 Now suppose $M$ is compact, $g$ is a K\"ahler metric, and $\phi:(M,g)\to (N, h)$ is a harmonic map. Then  \begin{align}\label{2.4}0&=\int_M \sqrt{-1}\, \p\pbar\varepsilon(\phi)\wedge\frac{\omega^{m-2}}{(m-2)!}\\\nonumber &=-\int_M  \sum_{\alpha,\beta=1}^m R^N(\phi_*Z_\alpha ,\, \phi_*Z_\beta,\, \phi_*\overbar{Z_\alpha},\, \phi_*\overbar{Z_\beta})
  \,d V_{\omega}+\int_M|D''\p \phi|^2 d V_{\omega}.\end{align}
 Hence if $(N,h)$ satisfies \eqref{1.1}, we must have $D''\p\phi=0$, i.e., $\phi$ is pluriharmonic. The preceding argument for Theorem \ref{theorem 1} parallels that in \cite[Proposition 3]{Siu}.\\ \indent The relation between the $(1,1)$-form $\varepsilon(\phi)$ and the energy $E(\phi)$ is as follows (see, e.g., \cite{Ohnita}):  
 Suppose $\phi: M\to N$ is a smooth map from an $n$-dimensional compact Hermitian manifold $(M,g)$ into a Riemannian manifold $(N,h)$. As before, we denote by $\omega$ the fundamental $2$-form of $(M,g)$. Then the energy of $\phi$ is given by\vspace{-0.05cm}\[E(\phi)=\int_M e(\phi)\, d V_\omega=\int_M \varepsilon(\phi)\wedge\frac{\omega^{n-1}}{(n-1)!},\]\vspace{-0.05cm}where $e(\phi)=\frac{1}{2}\,|d\phi|^2=\langle\varepsilon(\phi),\omega\rangle$. If $\phi$ is pluriharmonic, then \eqref{2.1} shows that $\varepsilon(\phi)$ is $d$-closed. Thus if, in addition, $g$ is a K\"ahler metric, then the number $E(\phi)$ depends only on the coho-mology classes $[\varepsilon(\phi)], [\omega]\in H^2(M,\RR)$. In particular, the energy of a pluriharmonic map $\phi:M\to N$ from a polarized compact K\"ahler manifold $(M, [\omega])$ into a Riemannian manifold $(N,h)$ is well defined.
 \smallskip
 \section{Main results} 
 The proof of Theorem \ref{theorem 2} is an application of the Bochner type identity \eqref{2.2} to noncompact K\"ahler manifolds. More specifically, we consider K\"ahler manifolds that admit a proper holomorphic submersion onto $\Delta$. We shall derive Theorem \ref{theorem 2} from Proposition \ref{proposition 7}, which gives a second variation formula  for the energy of pluriharmonic maps.
 \\\indent 
In Section \ref{section 3.1} we outline the proof of Theorem \ref{theorem 2}, in which we incorporate the main components of the proof of Proposition \ref{proposition 7}; this will shed light on the relation between the\linebreak arguments for Theorem \ref{theorem 1} and \ref{theorem 2}. In Section \ref{section 3.2} we extend a classical result of Lichnerowicz \cite{Lichnerowicz} (see Proposition \ref{proposition 4}), derive a closure property of holomorphic maps between K\"ahler mani-folds, and  proceed to establish Proposition \ref{proposition 7}. Finally, in Section \ref{section 3.3} we formulate sufficient conditions for the strict plurisubharmonicity of the energy.

\subsection{Outline of the proof of Theorem 2}\label{section 3.1}

\noindent{\bf Setting and reduction.}
First, by replacing $\Delta$ by a neighborhood of a point $t_0\in \Delta$ if necessary, we may assume that the polarization $\widetilde{\omega}$ of the family $\pi:\X\to \Delta$ is represented by a K\"ahler form $\omega_{\X}$ on $\X$ (see, e.g., \cite{Fujiki1, Fujiki2}). We choose local holomorphic coordinates on $\X$ such that $\pi$ is given in these coordinates by the projection $(z^1,...,z^n,t)\mapsto t$. For each $t\in\Delta$, the restriction of $\omega_{\X}$ to $M_t$ yields a K\"ahler form on $M_t$, which we denote by $\omega_t=\sqrt{-1}\,g_{\alpha\bbeta}(z^1,...,z^n,t)\,dz^\alpha\wedge d{\overbar{z^\beta}}$.

\noindent {\bf Fiber integral formula for $\boldsymbol{\sqrt{-1}\,\p\pbar E}$.} The smoothness of the family $\{\phi_t: M_t\to N\}_{t\in\Delta}$ means that the map $f:\X\to N$, $f(p)=\phi_{\pi(p)}(p),$ is smooth. Thus we can write the function $E:\Delta\to\RR$ as the fiber integral\[E=\int_{\X/\Delta}\varepsilon(f)\wedge\frac{\omega_{\X}^{n-1}}{(n-1)!},\] and hence \begin{equation}\label{3.1}
  \sqrt{-1}\,\p\overline{\p}E=\int_{\X/\Delta}\sqrt{-1}\,\p\overline{\p}\varepsilon(f)\wedge\frac{\omega_{\X}^{n-1}}{(n-1)!}.\end{equation}\vspace{-0.3cm}
  
  \noindent {\bf Computation of the fiber integral in (3.1).} Let $D$ be the exterior covariant derivative induced on $f^* T_{\CC} N$, and let $H$ be the horizontal lift of $\p/\p t$ with respect to $\omega_{\X}$.  We write $\lvert H \rvert$ for the pointwise norm of $H$ induced by $\omega_{\X}$, i.e., $|H|^2=-\sqrt{-1}\,\omega_\X(H,\overbar{H})$.  Applying \eqref{2.3}, we obtain 
 \begin{align}\label{3.2}
\sqrt{-1}\,\p\overline{\p}\varepsilon(f)\wedge\frac{\omega_{\X}^{n-1}}{(n-1)!}=\bigg\{&-2\, R^N\left( f_*\tfrac{\p\ }{\p z^\alpha},\,f_*H ,\, f_*\tfrac{\p\ }{\p \overbar{z^\beta}},\,f_*\overbar{H}\right)g^{\alpha\bbeta}\lvert H \rvert^{-2}\\\nonumber &+2\,h\left(D_{\overbar{H}}f_*\tfrac{\p\ }{\p z^{\alpha}},\, D_{H}f_*\tfrac{\p\ }{\p \overbar{z^\beta}}\right) g^{\alpha\bbeta}\lvert H \rvert^{-2}\bigg\}\,\frac{\omega_{\X}^{n+1}}{(n+1)!}.\end{align}
 \noindent From \eqref{3.1} and \eqref{3.2} it follows that \begin{align}\label{3.3}
    \frac{\p^2 E}{\p t\p \overline{t}}\ =\ &-2\int_{M_t} R^N\left(f_*\tfrac{\p\ }{\p z^\alpha}
    ,\,f_*H,\, f_*\tfrac{\p\ }{\p \overbar{z^\beta}},\,f_*\overbar{H} 
    \right) g^{\alpha\bbeta}\,d V_{\omega_t}\\\nonumber&+2\int_{M_t}h\left(  D_{\overbar{H}}f_*\tfrac{\p\ }{\p z^{\alpha}},\, D_{H}f_*\tfrac{\p\ }{\p \overbar{z^\beta}}\right) g^{\alpha\bbeta}\,d V_{\omega_t},\end{align}which is nonnegative on $\Delta$.

    \begin{remark} We summarize some of the main points that highlight the relation between the arguments for Theorem \ref{theorem 1} and \ref{theorem 2}. In Theorem \ref{theorem 1} the top degree form on $M$ given by \eqref{2.3} is posi-tive (in the sense of Lelong \cite{Lelong}), and 
     it follows from Stokes' theorem that $\phi$ is pluriharmonic. 
    In Theorem \ref{theorem 2} the top degree form on $\X$ given by \eqref{3.2} is also positive, and its fiber integral represents the Levi form of $E$ (cf.\ \eqref{2.4}),
    leading to the positivity of $\sqrt{-1}\, \p\pbar E$.\end{remark}
\subsection{Variations of the energy and fiber bundles}\label{section 3.2}

 Suppose $\phi: M\to N$ is a smooth map between K\"ahler manifolds, with $M$ compact. The differential $d\phi:T_\RR M\to \phi^*T_\RR N$ admits a decomposition $d\phi=d'\phi+d''\phi$ into complex linear and conjugate linear parts. In local holomorphic coordinates $(u^\mu)$ on $N$, we can write 
 $d'\phi=2\,\text{Re}\left(\p\phi^\mu\otimes\p/\p u^\mu\right)$ and $d''\phi=2\,\text{Re}\left(\pbar\phi^\mu\otimes\p/\p u^\mu\right)$, where $\text{Re}$ denotes the real part and $(\phi^\mu)$ is the coordinate representation of $\phi$ in $(u^\mu)$. Correspondingly, the preceding decomposition yields a decomposition of the energy $E(\phi)=E'(\phi)+E''(\phi)$. (Thus $\phi$ is holomorphic if and only if $E''(\phi)=0$.) It is a theorem of Lichnerowicz \cite{Lichnerowicz} that the number $K(\phi)=E'(\phi)-E''(\phi)$ is a smooth homotopy invariant. If $\{\phi_t:M\to N\}_{t\in(-1,1)}$ is a smooth family of harmonic maps, then $E(\phi_t)$ is constant; hence if, in addition, $\phi_0$ is holomorphic, then $\phi_t$ is holomorphic for every $t\in(-1,1)$.\\ 
 \indent  
Suppose now that $\pi:\X\to B$ is a smooth family of compact complex manifolds over the open unit ball $B$ in $\RR^k$, and that there exists a closed $2$-form $\omega_{\X}$ on $\X$ whose restriction to each fiber $M_t=\pi^{-1}(t)$ is a K\"ahler form on $M_t$, denoted by $\omega_t$. Let $\{\phi_t: M_t\to N\}_{t\in B}$ be a family of smooth maps into a fixed K\"ahler manifold $(N,h)$. As before, we assume that the family $\{\phi_t\}_{t\in B}$ is smooth, i.e., the map $f:\X\to N$, $f(p)=\phi_{\pi(p)}(p)$, is smooth. For each $t\in B$, we consider the numbers $K(\phi_t)$ and $E(\phi_t)$, defined with respect to the K\"ahler structures on\linebreak $M_t$ and $N$.   
\begin{proposition}\label{proposition 4}
 Under the preceding assumptions, $K(\phi_t)$ is constant. 
\end{proposition}
\vspace{-0.05cm}\noindent{\it Proof.} Let $n$ be the complex dimension of the fibers $M_t$, and let $\omega_N$ be the K\"ahler form on $N$ associated to $h$. Then \[K(\phi_t)=\int_{M_t}\phi_t^* \omega_N\wedge\frac{\omega_t^{n-1}}{(n-1)!}\]\\[-0.3cm] for every $t\in B$ (see, e.g., \cite{Eells, Lichnerowicz}). Since $\phi_t^* \omega_N=\left.\left(f^*\omega_N\right)\right|_{M_t}$ and $\omega_t=\omega_{\X}|_{M_t}$ for each $t\in B$,\linebreak the function $K: B\to \RR$, $t\mapsto K(\phi_t)$, is given by the fiber integral\[K=\int_{\X/B} f^*\omega_N\wedge\frac{\omega_{\X}^{n-1}}{(n-1)!}.\]
It follows that \[d K=\int_{\X/B}d\left(f^*\omega_N\wedge\frac{\omega_{\X}^{n-1}}{(n-1)!}\right)=0,\] since $\omega_N$ and $\omega_{\X}$ are closed.\qed\medskip

\indent Assume, in addition, that $\pi:\X\to B$ is a holomorphic family of compact complex mani-folds. In particular, $\X$ is a complex manifold, $k$ is even, $B$ is endowed with the standard complex structure, and $M_t$ is a complex submanifold of $\X$ for every $t\in B$.
\begin{proposition}\label{proposition 5}Under the above assumptions, the following hold:\begin{enumerate}[label=(\alph*)]
    \item If
 $f:\X\to N$ is pluriharmonic, then $E(\phi_t)$ is constant.\item If, in addition, $\phi_0$ is holomorphic, then $\phi_t$ is holomorphic for every $t\in B$.
\end{enumerate}
\end{proposition}
\vspace{-0.05cm}\noindent{\it Proof.} Let $n$ and $\omega_N$ be as in the proof of Proposition \ref{proposition 4}. For each $t\in B$, we have \[E(\phi_t)=\int_{M_t}\varepsilon(\phi_t)\wedge\frac{\omega_t^{n-1}}{(n-1)!}.\] Since $\varepsilon(\phi_t)=\varepsilon(f)\,|_{M_t}$ and $\omega_t=\omega_{\X}|_{M_t}$ for each $t\in B$, the function $E:B\to\RR$, $t\mapsto E(\phi_t)$, is given by the fiber integral\[E=\int_{\X/B}\varepsilon(f)\wedge\frac{\omega_{\X}^{n-1}}{(n-1)!},\] and we have \[dE=\int_{\X/B}d\varepsilon(f)\wedge\frac{\omega_{\X}^{n-1}}{(n-1)!}=0\]\\[-0.3cm] as $f$ is pluriharmonic. ($b$) holds since $2\,E''(\phi_t)= E(\phi_t)- K(\phi_t)$ is constant.\qed\pagebreak
\begin{remark} A similar analysis applies to a family $\{\phi_t:M_t\to N_t\}_{t\in B}$ of smooth maps between compact K\"ahler manifolds (in which case we assume that, for example, the family $\{N_t\}_{t\in B}$ is smooth and admits a closed $2$-form on its total space whose restriction to each fiber $N_t$ coincides with the K\"ahler form on $N_t$; see also (i) and (ii) in Proposition \ref{proposition 7}).\end{remark} 

We shall now complete the proof of Theorem \ref{theorem 2} by establishing the formulas stated in Section \ref{section 3.1}. The explicit expression \eqref{3.3} for $\p^2 E/\p t\p\overline{t}$ will be obtained as a special case of the following result:\vspace{0.05cm}
 \begin{proposition}\label{proposition 7}
 Let $\{\phi_t:M_t\to N\}_{t\in\Delta}$ be a family of smooth maps from compact K\"ahler manifolds $(M_t, g_t)$ into a fixed Riemannian manifold $(N,h)$. Suppose\vspace{-0.05cm} \begin{enumerate}[label=\normalfont(\roman*)]
     \item 
 $\{M_t\}_{t\in\Delta}$ is a holomorphic family of compact complex manifolds, with total space $\X$ and projection $\pi:\X\to\Delta$, \item there exists a K\"ahler form $\omega_{\X}$ on $\X$ such that $\omega_t=\omega_{\X}|_{M_t}$ is the K\"ahler form of $(M_t, g_t)$\linebreak for each $t\in \Delta$, and \item the family $\{\phi_t:M_t\to N\}_{t\in\Delta}$ is smooth.  \end{enumerate}\vspace{-0.05cm}
 Let $f:\X\to N$ be such that $f|_{M_t}=\phi_t$ for each $t\in\Delta$. If $\phi_0$ is pluriharmonic, then \begin{align}\label{3.4}
   \left. \frac{\p^2 E(\phi_t)}{\p t\p \overline{t}}\right|_{t=0}\ =\ &-2\int_{M_0} R^N\left(f_*\tfrac{\p\ }{\p z^\alpha}
    ,\,f_*H,\, f_*\tfrac{\p\ }{\p \overbar{z^\beta}}, \,f_*\overbar{H} 
    \right) g^{\alpha\bbeta}\,d V_{\omega_0}\\\nonumber&+2\int_{M_0}h\left(  D_{\overbar{H}}f_*\tfrac{\p\ }{\p z^{\alpha}},\, D_{H}f_*\tfrac{\p\ }{\p \overbar{z^\beta}}\right) g^{\alpha\bbeta}\,d V_{\omega_0}.\end{align}\end{proposition}
    
  \vspace{0.17cm}
   \begin{remark} $H$ is the horizontal lift of $\p/\p t$ with respect to $\omega_{\X}$, and $D$ is the exterior covariant derivative on $f^*T_\CC N$ induced by the Levi--Civita connection on $N$. The total space $\X$ is\linebreak covered by coordinate domains with local holomorphic coordinates such that $\pi$ is given in these coordinates by the projection $(z^1,...,z^n,t )\mapsto t$. In particular, for fixed $t$, $(z^1,...,z^n)$ are local holomorphic coordinates on $M_t$. For $t\in\Delta$, we write $\omega_t=\sqrt{-1}\,g_{\alpha\bbeta}(z^1,...,z^n,t)\,dz^\alpha\wedge d\overbar{z^\beta}$. 
   \end{remark}
   
  \begin{remark} Consider a harmonic map $\phi:M\to N$ between Riemannian manifolds $(M,g)$ and $(N,h)$, with $M$ closed and oriented. The formula for the second variation of the energy $E(\phi)$ is well known (see, e.g., \cite{Smith}). Proposition \ref{proposition 7} deals with the case where $\phi:M\to N$ is a pluriharmonic map from a compact K\"ahler manifold $(M, g)$, and formula \eqref{3.4} incorporates the effect of deforming $(M,g)$ in the sense of (i) and (ii). 
  \end{remark} 
   
   \vspace{-0.05cm}\noindent{\it Proof.} As in the proof of Proposition \ref{proposition 5},  we can write 
   the function $E:\Delta\to\RR$, $t\mapsto E(\phi_t)$, as\linebreak the fiber integral\vspace{-0.05cm}\[E=\int_{\X/\Delta}\varepsilon(f)\wedge\frac{\omega_{\X}^{n-1}}{(n-1)!}.\]\\[-0.36cm] Since $\pi:\X\to \Delta$ is holomorphic, it follows that\vspace{-0.05cm} \begin{equation}\label{3.5}
   \sqrt{-1}\,\p\pbar E=\int_{\X/\Delta}\sqrt{-1}\,\p\pbar\varepsilon(f)\wedge\frac{\omega_{\X}^{n-1}}{(n-1)!}.
   \end{equation}  
   \indent  
 For any $t\in\Delta$, if $\phi_t: M_t\to N$ is pluriharmonic, then at each point of $M_t$ we have  \begin{equation}\label{3.6}R^N\left(f_*\tfrac{\p\ }{\p z^\alpha},\,f_*\tfrac{\p\ }{\p z^\gamma}, f_*\tfrac{\p\ }{\p \overbar{z^\beta}},\,f_*\tfrac{\p\ }{\p \overbar{z^\delta}}\right)=0
 \ \ \ \ \text{and}\ \ \ \  
 D_{\frac{\p\ }{\p \overbar{z^\beta}}}f_*\tfrac{\p\ }{\p z^\alpha}=\left(\iota_t^*D\right)_{\frac{\p\ }{\p \overbar{z^\beta}}}(\phi_t)_*\tfrac{\p\ }{\p z^\alpha}=0, \end{equation}  where $\iota_t: M_t\hookrightarrow\X$ is the inclusion, for all $\alpha, \beta, \gamma, \delta\in\{1,...,n\}$ (see \cite[Lemma 1.2]{Ohnita}).\newpage\indent The holomorphic tangent bundle $T_{\X}$ admits an orthogonal decomposition $T_{\X}=$ $T_V\oplus T_H$ induced by $\omega_{\X}$, 
  where $T_V$ is the kernel of $\pi_*$. In particular, $\left. T_V\right|_{M_t}= T_{M_t}$ for each $t\in\Delta$, $H\in C^\infty(\X, T_H)$ is such that $\pi_*H=\p/\p t$, and $T_\X$ has a local frame given by $\big\{\frac{\p}{\p z^1},...,\frac{\p}{\p z^n}, H\big\}$.
  \\\indent 
 We shall now apply \eqref{2.3} with $\phi= f$ and $\omega=\omega_\X$. As before, $D'$ and $D''$ will denote the $(1,0)$- and $(0,1)$-components of $D$, respectively. We have, for example, \[
     \left(D''\p f\right)\left( \overbar{H},\tfrac{\p\ }{\p z^\alpha}\right) =  \left(D''_{\overbar{H}}\hspace{0.03cm} \p f\right)\left(\tfrac{\p\ }{\p z^\alpha}\right)  = D''_{\overbar{H}} \left(\p f\left(\tfrac{\p\ }{\p z^\alpha}\right)\right)-\p f\left(\pbar_{\overbar{H}} \tfrac{\p\ }{\p z^\alpha}\right) = D_{\overbar{H}}  f_*\tfrac{\p\ }{\p z^\alpha} \] and  \[\left(D''\p f\right)\left(\tfrac{\p\ }{\p \overbar{z^\beta}}, H\right) =  \left(D'\pbar f\right)\left(H,\tfrac{\p\ }{\p \overbar{z^\beta}}\right)\ =\ \overbar{D_{\overbar{H}}  f_*\tfrac{\p\ }{\p z^\beta}}= D_H  f_*\tfrac{\p\ }{\p \overbar{z^\beta}}\]\\[-0.38cm]
     for all $\alpha,\beta\in\{1,...,n\}$. Applying \eqref{2.3}, we obtain\begin{equation}\label{3.7}
     \sqrt{-1}\,\p\pbar\varepsilon(f)\wedge\frac{\omega_{\X}^{n-1}}{(n-1)!}= F\cdot \frac{\omega_{\X}^{n+1}}{(n+1)!},\end{equation} 
      where $F$ is a smooth real-valued function on $\X$ that satisfies
      \begin{align}\label{3.8}
      F|_{M_0}\ = \ &-2\,R^N\left(f_*\tfrac{\p\ }{\p z^\alpha},\, f_*H,\, f_*\tfrac{\p\ }{\p \overbar{z^\beta}},\, f_*\overbar{H}\right) g^{\alpha\bbeta}|H|^{-2}\\\nonumber&+2\,h\left(D_{\overbar{H}}f_*\tfrac{\p\ }{\p z^\alpha},\,D_{H} f_*\tfrac{\p\ }{\p \overbar{z^\beta}}\right)g^{\alpha\bbeta}|H|^{-2}. \end{align} 
      Indeed, the full expression of $F$ is given by
      \begin{align*}
      F\ =\ &-R^N\left(f_*\tfrac{\p\ }{\p z^\alpha},\, f_*\tfrac{\p\ }{\p z^\gamma},\, f_*\tfrac{\p\ }{\p \overbar{z^\beta}},\, f_*\tfrac{\p\ }{\p \overbar{z^\delta}}\right) g^{\alpha\bbeta}g^{\gamma\overline{\delta}}\\&-2\,R^N\left(f_*\tfrac{\p\ }{\p z^\alpha},\, f_*H,\, f_*\tfrac{\p\ }{\p \overbar{z^\beta}},\, f_*\overbar{H}\right) g^{\alpha\bbeta}|H|^{-2}\\&+h\left(D_{\frac{\p\ }{\p \overbar{z^\delta}}}f_*\tfrac{\p\ }{\p z^\alpha},\,D_{\frac{\p\ }{\p z^\gamma}} f_*\tfrac{\p\ }{\p \overbar{z^\beta}}\right)g^{\alpha\bbeta}g^{\gamma\overline{\delta}}\\&+2\,h\left(D_{\overbar{H}}f_*\tfrac{\p\ }{\p z^\alpha},\,D_{H} f_*\tfrac{\p\ }{\p \overbar{z^\beta}}\right)g^{\alpha\bbeta}|H|^{-2}\\&+\left\vert\left(D''\p f\right)\left(\overbar{H}, H\right)\right\vert^2_h\, |H|^{-4}\\&-\Big\vert\, g^{\alpha\bbeta}D_{\frac{\p\ }{\p \overbar{z^\beta}}}f_*\tfrac{\p\ }{\p z^\alpha}+|H|^{-2}\left(D''\p f\right)\left(\overbar{H}, H\right)\Big\vert^2_h,\end{align*} and at each point of $M_0$, \eqref{3.6} implies that the first and the third terms on the right side vanish and that the last two terms cancel out. Thus we obtain \eqref{3.8}.\\\indent
       Write $\psi=|H|^2$. We have \begin{equation}\label{3.9} \frac{\omega_{\X}^{n+1}}{(n+1)!}=\sqrt{-1}\,\psi\, dt\wedge d\overline{t}\wedge\frac{\omega_{\X}^{n}}{n!}.\end{equation} This can be verified directly: If $\{\eta^1,...,\eta^n, dt\}$ is the local coframe dual to $\big\{\frac{\p}{\p z^1},...,\frac{\p}{\p z^n}, H\big\}$, then \[\omega_{\X}=\sqrt{-1}\,g_{\alpha\bbeta}\,\eta^\alpha\wedge\overbar{\eta^\beta}+\sqrt{-1}\,\psi\,dt\wedge d\overline{t},\] and both sides of \eqref{3.9} are equal to \[\left(\sqrt{-1}\right)^{n+1}\psi\det\big(g_{\alpha\bbeta}\big)\,dt\wedge d\overline{t}\wedge\eta^1\wedge \overbar{\eta^1}\wedge...\wedge\eta^n\wedge \overbar{\eta^n}.\]
      \indent Combining \eqref{3.5}, \eqref{3.7}, and \eqref{3.9}, we obtain  \begin{align*}
    \sqrt{-1}\,\p\pbar E\ &=\ \int_{\X/\Delta}\sqrt{-1}\,\p\pbar \varepsilon(f)\wedge\frac{\omega_{\X}^{n-1}}{(n-1)!}\ =\ \int_{\X/\Delta}F\cdot\frac{\omega_{\X}^{n+1}}{(n+1)!}\\&=\ \int_{\X/\Delta}\sqrt{-1}\,F\,\psi\, dt\wedge d\overline{t}\wedge\frac{\omega_{\X}^n}{n!}.\end{align*}Since $dt=\pi^*(dt)$, it follows that \[
    \sqrt{-1}\,\p\pbar E (0)\ =\ \sqrt{-1}\, dt\wedge d\overline{t}\cdot\int_{M_0}\left.\left(F\,\psi\,\frac{\omega_{\X}^{n}}{n!}\right)\right|_{M_0}.\] From \eqref{3.8} we see that $\left.\left(F\,\psi\,\frac{\omega_{\X}^{n}}{n!}\right)\right|_{M_0}$ is the top degree form on $M_0$ given by 
    \begin{align*}
     &-2\, R^N\left(f_*\tfrac{\p\ }{\p z^\alpha},\, f_*H,\, f_*\tfrac{\p\ }{\p \overbar{z^\beta}},\, f_*\overbar{H}\right) g^{\alpha\bbeta}\, d V_{\omega_0}\\&+2\,h\left(D_{\overbar{H}}f_*\tfrac{\p\ }{\p z^\alpha},\,D_{H} f_*\tfrac{\p\ }{\p \overbar{z^\beta}}\right)g^{\alpha\bbeta}\, dV_{\omega_0}. \end{align*}Thus we conclude that \eqref{3.4} holds.
     \qed\medskip
     
    \subsection{Strict plurisubharmonicity}\label{section 3.3}
    
    In this section, we provide sufficient conditions for the strict plurisubharmonicity of the energy. We use the notation of Proposition \ref{proposition 7}, and  consider the function $E:\Delta\to\RR, t\mapsto E(\phi_t)$. 
   \begin{corollary}\label{corollary 10}
   If $(N,h)$ satisfies the curvature condition \eqref{1.1}, and if \begin{equation}\label{3.10}
    \left.\left(d\overbar{z^\beta}\otimes D_H f_*\tfrac{\p\ }{\p \overbar{z^\beta}}\right)\right|_{M_0}\end{equation} is not identically zero, then $E$ is strictly subharmonic at $t=0$.\end{corollary}\vspace{-0.05cm}\noindent {\it Proof.} The assumptions show that the first term on the right side of \eqref{3.4} is nonnegative, and that the second term is  positive,  the integrand in the second term being the pointwise norm of  \eqref{3.10} induced by $g_0$ and $h$.\qed\medskip
    
   \indent We now proceed to examine \eqref{3.10}. For simplicity we write $M=M_0$ and $\phi=\phi_0$. Let \[\rho:T_{\Delta,\hspace{0.02cm}0}\to H^1(M,\mathcal{O}(T_{M}))\] be the Kodaira--Spencer map of $\pi$ at $t=0$. Write $H={\p}/{\p t}+a^\alpha{\p}/{\p z^\alpha}$, let $A_{\bbeta}^\alpha={\p a^\alpha}/{\p \overbar{z^\beta}}$ for $\alpha, \beta\in\{1,...,n\}$, and let \[A=\left.\left(\pbar    H\right)\right|_{M}=A_{\bbeta}^\alpha\,d \overbar{z^\beta}\otimes\tfrac{\p\ }{\p z^\alpha}.\]  The cohomology class $[A]$ in $ H^{0,1}(M, T_{M})$ corresponds to the Kodaira--Spencer class $\rho ({\p}/{\p t})$ under the Dolbeault isomorphism. We shall use the notation  \[\phi_*A = A_{\bbeta}^\alpha\,d \overbar{z^\beta}\otimes d\phi\left(\tfrac{\p\ }{\p z^\alpha}\right).\]\pagebreak Let $\nabla''$ be the $(0,1)$-component of the exterior covariant derivative induced on $\phi^*T_\CC N$. Since \[D_{\frac{\p\ }{\p \overbar{z^\beta}}}f_*H- D_{H}f_*\tfrac{\p\ }{\p \overbar{z^\beta}}=f_*\left[\tfrac{\p\ }{\p \overbar{z^\beta}}, H\right]= A_{\bbeta}^\alpha \,f_*\tfrac{\p\ }{\p z^\alpha}\] for all $\beta\in\{1,...,n\}$, we have \begin{align*}
     \left.\left(d\overbar{z^\beta}\otimes D_H f_*\tfrac{\p\ }{\p \overbar{z^\beta}}\right)\right|_{M}&=\iota^*\left(d\overbar{z^\beta}\otimes D_{\frac{\p\ }{\p \overbar{z^\beta}}}f_*H-d\overbar{z^\beta}\otimes A_{\bbeta}^\alpha \,f_*\tfrac{\p\ }{\p z^\alpha}\right)\\&=\nabla''(\iota^*f_*H)-\phi_* A, \end{align*}where $\iota:M\hookrightarrow \X$ is the inclusion. Thus Corollary \ref{corollary 10} can be restated as follows.
      
     \begin{corollary}\label{corollary 11} If $(N, h)$ satisfies the curvature condition \eqref{1.1}, and if $\phi_* A\neq \nabla''(\iota^*f_*H)$, then $E$ is strictly subharmonic at $t=0$. \end{corollary}
     
    \indent 
    We shall consider situations where the operator $\nabla''$ satisfies  $(\nabla'')^2=0$, and hence induces the structure of a holomorphic vector bundle on $\phi^*T_\CC N$. The pluriharmonic map equation $\nabla''\p\phi=0$ then shows that $\p\phi$ is a holomorphic section of the bundle $ T^*_{M}\otimes\phi^*T_\CC N$. Since \[(\nabla'')^2\,(V,W)=R^N(d\phi(V), d\phi(W))\ \ \ \ \text{and}\ \ \ \ \text{dim}_\CC d\phi(\overbar{T_{M,p}})=\text{dim}_\CC d\phi({T_{M,p}})\]  for every $V,W\in \overbar{T_{M,p}}$ and for every $p\in M$, we see that the integrability condition $(\nabla'')^2=0$ holds whenever $\text{dim}_\CC d\phi({T_{M,p}})\leq 1$ for every $p\in M$. 
     \begin{corollary}\label{corollary 12} If the curvature operator of $(N, h)$ is nonpositive, and if the cohomology class $[\phi_*A]$ in $H^{0,1}(M,\phi^*T_\CC N)$ is not zero, then $E$ is strictly subharmonic at $t=0$.
       \end{corollary} 
      \vspace{-0.05cm}\noindent{\it Proof.} This is an immediate consequence of \cite[Proposition 1.4]{Ohnita} and Corollary \ref{corollary 11}.
      \qed\medskip\vspace{0.03cm}
      
       Motivated by  \cite[Theorem 3]{Toledo}, we establish a variant of the preceding corollary.
     \begin{corollary}\label{corollary 13} Suppose that $(N, h)$ satisfies the curvature condition: \begin{equation}\label{3.11}
     R^N(X,Y,\overline{X}, \overline{Y})< 0 \ \ \, \text{for  all $X,Y\in T_{\CC,p }N$ with $X\wedge Y\neq 0$ and for all $p\in N$.}    
     \end{equation} Suppose $d\phi$ is nowhere zero on $M$. If the cohomology class $[\phi_*A]$ in $H^{0,1}(M,d\phi(T_{M}))$ is not zero, then $E$ is strictly subharmonic at $t=0$.\\[-0.38cm]
       \end{corollary}
       
       \begin{remark} 
         By \cite[Lemma 1.2]{Ohnita}, \eqref{3.11} implies that $\text{dim}_\CC d\phi({T_{M,p}})\leq 1$ for all $p\in M$. Since $d\phi$  is nowhere zero on $M$, we have   $\text{dim}_\CC d\phi({T_{M,p}})= 1$ for all $p\in M$.  It follows that $d\phi(T_M)$, being the image of $\p\phi$, is a holomorphic line subbundle of $\phi^* T_\CC N$. 
       \end{remark}
 \vspace{-0.05cm}\noindent{\it Proof.} Assume, to reach a contradiction, that $\sqrt{-1}\,\p\pbar E(0)=0$. Then \eqref{3.4} and \eqref{3.11} show that  at each point of $M$,  \[f_*\tfrac{\p\ }{\p z^\alpha}\wedge f_*H=0\] for all $\alpha\in\{1,...,n\}$; equivalently, 
 \[\iota^* f_*H\in C^\infty(M, d\phi(T_M)).\] On the other hand, Corollary \ref{corollary 11} implies that $\phi_* A= \nabla''(\iota^*f_*H)$. Thus the cohomology class $[\phi_*A]$ in $H^{0,1}(M,d\phi(T_{M}))$ is zero, and we have a contradiction.
 \qed
   
\section*{Acknowledgements} The author is indebted to his advisor Sai-Kee Yeung for stimulating conversations during the writing of this paper, and to Laszlo Lempert for fruitful discussions on plurisubharmonic functions and for research funding support during the summer of 2021.

\end{document}